\documentclass[12pt]{article}
\usepackage{caption}
\usepackage{graphicx,latexsym,amssymb,amsmath}
\usepackage{tikz}
\usepackage{tabularx}
\usepackage{diagbox}
\definecolor{myblue}{RGB}{173, 216, 230}
\DeclareMathOperator*{\argmin}{arg\,min}
\usepackage[french,english]{babel}
\usepackage{enumerate}
\usepackage{here} 
\def\R{{\mathbb{R}}}

\def\P{{\Pr}}
\def\E{{\mathbb{E}}}

\newcommand{{\convp}}{{\buildrel p\over\longrightarrow}}

\newcommand{\red}[1]{{\color{black} #1}}

\newcommand{{\Vs}}{{\cal V}}
\newcommand{{\Ps}}{{\cal P}}
\newcommand{{\Ss}}{{\cal S}}
\newcommand{{\Xs}}{{\cal X}}
\newcommand{{\Ls}}{{\cal L}}
\newcommand{{\Ns}}{{\cal N}}
\newcommand{{\Zs}}{{\cal Z}}
\newcommand{{\Fs}}{{\cal F}}

\newtheorem{Assumption}{Assumption}[section]
\newtheorem{Theorem}{Theorem}[section] 
\newcommand{\proofend}{$\quad\Box{~}$}

\renewcommand{\baselinestretch}{1.1}

\makeatletter

\@addtoreset{equation}{section}
\makeatother

\usepackage{natbib}

\begin{document}
\title{\bf Two-stage least squares with a randomly right censored outcome}

\author{
{\large Jad B\textsc{eyhum}}
\footnote{ORSTAT, KU Leuven. Naamsestraat 69, 3000 Leuven (Belgium).\newline\indent  Financial support from the European Research Council (2016-2021, Horizon 2020 / ERC grant agreement\newline \indent  No.\ 694409) is gratefully acknowledged.\newline \indent Declarations of interest: none.}\\\texttt{\small jad.beyhum@kuleuven.be}
}

\date{\today}

\maketitle

\begin{abstract}

This note develops a simple two-stage least squares (2SLS) procedure to estimate the causal effect of some endogenous regressors on a randomly right censored outcome in the linear model. The proposal replaces the usual  ordinary least squares regressions of the standard 2SLS by weighted least squares regressions. The weights correspond to the inverse probability of censoring. We show consistency and asymptotic normality of the estimator. The estimator exhibits good finite sample performances in simulations. 

\end{abstract}
\smallskip


\smallskip 

\bigskip

\def\baselinestretch{1.3}

\newpage
\normalsize

\setcounter{footnote}{0}
\setcounter{equation}{0}
\section{Introduction}
Two-stage least squares (2SLS) is arguably the most famous estimator using instrumental variables (see \cite{wooldridge2010econometric}). It can be used to treat omitted variable bias when instrumental variables are available. The main reason why 2SLS is popular among practitioners is its simplicity. It requires a finite number of ordinary least squares regressions, which ensures that it possesses advantegeous computational and finite sample properties. In the contexts where 2SLS would be useful, there are many occurrences where the outcome of interest is right censored. For instance, unemployment durations or the duration until recovery from a disease are often not observed after the end of the study. It is therefore important to develop a 2SLS procedure for this type of data.

This paper develops an analogous 2SLS procedure for right censored outcomes. We propose to replace the usual ordinary least squares regressions of 2SLS with weighted least squares regressions. For noncensored observations, the weights are the inverse of the probability of censoring estimated by the Kaplan-Meier estimator. For censored observations, the weights are equal to $0$. The method retains the simplicity of the usual 2SLS.  It is assumed that the censoring is independent of the other variables of the model and that there is sufficient follow-up. We show the consistency and the asymptotic normality of the estimator. The estimator exhibits excellent finite sample performances in simulations. 

There exists an extensive literature on instrumental variable methods with randomly right censored duration outcomes. \cite{BR}, \cite{tchetgen2015instrumental}, \cite{li2015instrumental}, \cite{chan2016reader} and \cite{chernozhukov2015quantile} estimate causal effects in semiparametric frameworks, which do not correspond to the linear model of our paper. \cite{frandsen2015treatment},  \cite{richardson2017nonparametric}, \cite{blanco2019bounds}, \cite{beyhum2021nonparametric}, \cite{sant2016program}, \cite{sant2020nonparametric},  and \cite{centorrino2021nonparametric} provide nonparametric estimation results. The last three papers of the latter list are especially close to the present work. They also use inverse probability of censoring weighting. However, none of these papers considers the simple linear model as the present work does. The advantage of the proposed procedure over the aforementioned articles is therefore its simplicity, which may make it more popular with practitioners. Finally, like the present paper,  \cite{kjaersgaard2016instrumental} considers parametric models. Nevertheless, they use a pseudo-observation approach and only show asymptotic unbiasedness of the estimator. In contrast, our proposed procedure is simpler than the pseudo-observation approach and we show asymptotic normality. 
\section{The model}
We are interested in the effect of a $K$-dimensional random vector of regressors $X$ on a real-valued outcome $T$. The two variables are related by the following model
\begin{equation}\label{model}T=X^\top \beta_0+U,
\end{equation}
where $\beta_0\in\R^K$ is the parameter of interest and $U$ is a real-valued random error term. The variable $X$ is endogenous but we possess an instrumental variable $Z$ with support in $\R^L$. It satisfies the following conditions, which ensure consistency of the 2SLS estimator in the absence of censoring.
\begin{Assumption}\label{IV} The following holds:
\begin{itemize}
\item[(i)] $E[||Z||^2]+E[||X||^2]+ E[||Y||^2]<\infty$;
\item[(ii)] $\E[ZU]=0$;
\item[(iii)] $\E[ZZ^\top]$ is positive definite;
\item[(iv)] $\textrm{rank}(\E[ZX^\top])=K$.
\end{itemize}
\end{Assumption}
By $||\cdot||$, we denote the euclidian norm.
Unlike in the standard linear IV model, the outcome variable $T$ is right censored by a real-valued censoring variable $C$. Let $Y=\min(T,C)$ and $\delta=I(T\le C)$. The observables are $(Y,\delta,X,Z)$. In a labor economics example, $T$ could be the unemployment, duration, $C$ the time between the beginning of the unemployment spell and the end of the study, $X$ a treatment indicator (e.g. participation in a job training) and $Z$ a randomized treatment/control group assignment indicator.

\section{Estimation procedure}
\subsection{Stute's estimator for moments}
For estimation, we assume that we possess an i.i.d. sample $\{Y_i,\delta_i,X_i, Z_i\}_{i=1}^n$ drawn from the distribution of $(Y,\delta,X, Z)$. Because the outcome is right censored, the second-step of the usual two-stage least squares estimator is not feasible. \cite{Stute1993consistent} proposes a generic method to estimate moments when one variable is censored. Let $\varphi: \R^K\times \R^L\times \R \mapsto \R^M$ and suppose that the goal is to estimate $\E[\varphi(X,Z,T)]$. The usual estimator 
$ n^{-1}\sum_{i=1}^n \varphi(X_i,Z_i,T_i)$ cannot be computed because $T_i$ is censored for some observations. The idea of \cite{Stute1993consistent} is to use the following equality 
\begin{equation}\label{stute_moment}  \E[\varphi(X,Z,T)]=\E\left[\frac{\delta}{G(Y)}\varphi(X,Z,Y)\right],\end{equation}
where $G(t)=\P(C\le t)$ is the cumulative distribution function of $C$.
It says that the expectation of $\varphi(X,Z,T)$ is equal to the expectation of $\varphi(X,Z,Y)$ weighted by $\delta/G(Y)$. Under standard assumptions from the survival analysis literature, $G$ is identified and therefore all elements of the right-hand-side of \eqref{stute_moment} are observed, so that this equation can be used to construct a consistent estimator. 

Let us now give a formal definition of \cite{Stute1993consistent} estimator. We assume that the distribution of $Y$ is continuous, which avoids ties between observations and simplifies the definition of Stute's estimator. Some notations are needed. $Y_{(1)}\le\dots\le Y_{(n)}$ denote the order statistics of $\{Y_i\}_{i=1}^n$. For a vector $R\in\R^n$, $R_{(i)}$ is the coordinate of $R$ associated with $Y_{(i)}$.  Let us introduce the Kaplan-Meier weights:
\begin{equation} \label{kmweights} w_{(1)}= \frac{\delta_{(1)}}{n},\ w_{(i)}= \frac{\delta_{(i)}}{n-i+1}\prod_{j=1}^{i-1}\left( \frac{n-j}{n-j+1}\right)^{\delta_{(j)}},\ i=1,\dots,n,
\end{equation}
It can be shown that $w_{(i)}=\delta_{(i)}/\widehat{G}(Y_{(i)})$, where $\widehat{G}$ is the Kaplan-Meier estimator of $G$. Hence, for uncensored observations, the weights correspond to an estimate of the inverse of the probability that an observation is censored. Stute's estimator is 
\begin{equation}\label{stute_est} \widehat\E_w[\varphi(X,Z,Y)]= \sum_{i=1}^n w_{(i)}\varphi(X_{(i)},Z_{(i)}, Y_{(i)}).\end{equation}
\cite{Stute1993consistent} shows that $\widehat\E_w[\varphi(X,Z,Y)]$ converges in probability to $\E[\varphi(X,Z,T)]$. \cite{Stute1993consistent} gives conditions under which $\sqrt{n}(\widehat\E_w[\varphi(X,Z,Y)]-\E[\varphi(X,Z,T)])$ is asymptotically normal.

\subsection{Estimation procedure}
The two steps of 2SLS can be expressed as the minimization of some moments, therefore it is natural to use Stute's estimator to estimate these moments. The proposed estimation procedure is as follows. The first step corresponds to the weighted by $\{w_{(i)}\}$ regression of $X $ on $Z$, that is
\begin{equation}\label{estimator1step}\widehat{\Gamma}= \argmin_{\Gamma \in \R^{L\times K}}\sum_{i=1}^n w_{(i)}||X_{(i)}-Z_{(i)}\Gamma||^2.\end{equation}
Two remarks are in order. First, by standard arguments 
\begin{equation}\label{estimator1step2}\widehat{\Gamma}=\left(\sum_{i=1}^n w_{(i)}Z_{(i)}Z_{(i)}^\top\right)^{-1}\left(\sum_{i=1}^n w_{(i)}Z_{(i)}X_{(i)}^\top\right).\end{equation}
This expression is useful to compute $\widehat{\Gamma}$ in practice and to derive its asymptotic properties. Second, since neither $X$ nor $Z$ are censored, we could in principle use the ordinary least squares regression of $X$ on $Z$ as the first step. However, although this alternative may seem more natural, it complexifies the asymptotic theory of the estimator since different weights would be used in the first and second-step estimator.

The second-step estimator is
\begin{equation}
\label{est}\widehat{\beta}  =\argmin_{\beta\in\R^K} \sum_{i=1}^n w_{(i)}\left[Y_{(i)}-(Z_{(i)}^\top \widehat{\Gamma})\beta\right]^2.
\end{equation}
Notice that we have 
\begin{equation}\label{closed-form}
\widehat{\beta}= \left[\widehat{\Gamma}^\top\left(\sum_{i=1}^n  w_{(i)}Z_{(i)}Z_{(i)}^\top\right)\widehat{\Gamma}\right]^{-1} \widehat{\Gamma}^\top\left(\sum_{i=1}^n w_{(i)}Z_{(i)}Y_{(i)}\right)\end{equation}
which allows to compute $\widehat \beta$ easily.

\section{Asymptotic properties}

\subsection{Consistency} For a random variable $R$, let $\tau_R=\sup\{t\in\R:\ \P(T\le t)<1\}.$ The quantity $\tau_R$ is the upper bound of the support of $R$.
We make the following assumptions on the censoring variable
\begin{Assumption}\label{censoring} The following holds:
\begin{enumerate}
 \item[(i)] The censoring variable $C$ is independent of $(T,X)$;
\item[(ii)] $\tau_{T}< \tau_{C}$.
 \end{enumerate}

 \end{Assumption}
  The first condition is the usual independent censoring assumption from the survival analysis literature. Condition (ii) is a sufficient follow-up assumption, it means that the largest possible censoring time is larger than the largest possible duration $T$.  Notice that it implies that the duration is bounded from above. Assumption \ref{censoring} can be relaxed, for instance  \cite{Stute1993consistent} and \cite{stute1996distributional} impose weaker but less understandable conditions. 
 We have the following consistency result. 
 \begin{Theorem}\label{const}Under Assumptions \ref{IV} and \ref{censoring}, $\widehat{\beta}$ converges in probability to $\beta_0$. 
 \end{Theorem}

 \subsection{Asymptotic normality}   Recall that $G(t)=\P(C\le t)$. Moreover, for $t\in\R_+$, let $G(t-)$ be the left limit of $G$ at $t$ and $H(t)=\P(Y\le t)$ be the cumulative distribution function of $Y$ .  Next, for $\ell=1,\dots, L$, let 
 $$\varphi_\ell:\begin{array}{cl}\R^{K}\times \R^L \times \R& \mapsto \R \\
  (x,z,t) &\mapsto z_{\ell}(t-x^\top\beta_0).\end{array}$$ To ensure the asymptotic normality of the estimator, we impose the following condition.
This is a condition on the tail of the censoring distribution. It imposes that the tail cannot be too heavy.
\begin{Assumption}\label{stute_rate} It holds that
\begin{align}
\label{cond_stute_1}
 \int \int \int_0^{\tau_{T}}  |\varphi_k(x,z,s)|\sqrt{\int_0^{s-} \frac{dG(y)}{(1-H(y))(1-G(y-))}} dF(x,z, s)<\infty,
\end{align}
 Our estimator has a complicated asymptotic variance. In order to define it, we introduce the following quantities. For $x\in\R^K$, $z\in\R^{L}$ and $t\in\R$, let 
 \begin{align*}
 F(x, z,  t) &= \P( X_{1}\le x_1,\dots, X_{K}\le x_k,Z_1\le z_1,\dots,Z_L\le z_L,T\le  t)\\
 H^0(y)&= \P(T\le t,\delta=0)\\
 H^{11}(x,e,t)&=\P( X_{1}\le x_1,\dots, X_{K}\le x_k,Z_1\le z_1,\dots,Z_L\le z_L,T\le  t,\delta=1).
 \end{align*} 

for all $k=1,\dots, p$.
\end{Assumption}
For a function $\varphi:\R^K\times \R^L \times \R\mapsto \R$ and $t\in \R$, we define
\begin{align*}
\gamma_1^\varphi(t)&=\frac{1}{1-H(t)}\int \int \int 1_{\{t<s\}} \varphi(x,z,s) \frac{dH^{11}(x,z,s)}{1-G(s-)};\\
\gamma_2^\varphi(t)&=\int \int \int  \int 1_{\{s<t, s<y\}} \varphi(x,z,s) \frac{dH^0(s)dH^{11}(x,z,y)}{(1-H(s))^2(1-G(y-))}.
\end{align*}
For $\ell=1,\dots, L$, we define $\psi_\ell= \varphi_\ell(X, Z, Y)(1-G(Y-))^{-1}\delta+ \gamma^{\varphi_\ell}(Y)(1-\delta) -\gamma_2^{\varphi_\ell}(Y)$ and $\Psi=(\psi_1,\dots,\psi_L)^\top$.  The following result states that our estimator is asymptotically normal.
 \begin{Theorem}\label{AN} Let the assumptions of Theorem \ref{const} and Assumption \ref{stute_rate} hold. Then, we have 
 $$ \sqrt{n}\left(\widehat \beta -\beta_0\right)\xrightarrow{d}\mathcal{N}(0,\Sigma),$$
 where $\Sigma =W\Sigma_{\Psi}W^\top $ with 
 \begin{align*}
 W&= \left[ \E[ZX^\top] \E[ZZ^\top] ^{-1}  \E[ZX^\top]\right]^{-1}\E[ZX^\top]\E[ZZ^\top]^{-1};\\
\Sigma_{\Psi} &= E[\Psi\Psi^\top].
 \end{align*}
 \end{Theorem}
 
 \subsection{Estimation of the asymptotic variance} 
 Let us propose an estimator of the asymptotic variance of the estimator. Its accuracy is demonstrated in simulations in the next section. A natural estimator for $W$ is 
 $$\widehat{W}=  \left[\widehat{\Gamma}^\top\left(\sum_{i=1}^n  w_{(i)}Z_{(i)}Z_{(i)}^\top\right)\widehat{\Gamma}\right]^{-1} \widehat{\Gamma}^\top.$$
 To estimate $\Sigma_{\Psi}$, for $i=1,\dots, n$, we set $\widehat{U}_{(i)} = Y_{(i)} -X_{(i)}^\top\widehat{\beta}$. Then, for $\ell=1,\dots,K$, we let
\begin{align*}
\widehat{\gamma}_1^{\varphi_\ell}(t)&=\frac{1}{1-\widehat{H}(t)}\frac{1}{n}\sum_{i=1}^n 1_{\{t<Y_{(i)} \}}  \frac{ \delta_{(i)}Z_{(i)\ell}\widehat{U}_{(i)}}{1-\widehat{G}(Y_{(i)}-)};\\
\widehat{\gamma}_2^{\varphi_\ell}(t)&= \red{\frac{1}{n^2}}\sum_{i=1}^n \sum_{j=1}^n \int 1_{\{Y_{(j)}<t, Y_{(j)}<Y_{(i)}\}}   \frac{(1-\delta_{(j)})\delta_{(i)}Z_{(i)\ell}\widehat{U}_{(i)}}{(1-\widehat{H}(Y_{(j)}))^2(1-\widehat{G}(Y_{(i)}-))},
\end{align*}
where $\widehat{H}(t)=n^{-1}\sum_{i=1}^n 1_{\{Y_i\le t\}}$ and$\widehat{G}$ is the Kaplan-Meier estimator of $G$, that is 
$$ G(t) = 1-\prod_{i\in\{1,\dots,n\}:\ Y_{(i)}\le t}\left(\frac{1}{n-i+1}\right)^{1-\delta_{(i)}}$$
Also, we write $$\widehat{\psi}_{i\ell}= Z_{(i)\ell}\widehat{U}_{(i)}\delta_{(i)}(1-\widehat{G}(Y_{(i)}-))^{-1}+\widehat{ \gamma}_1^{\varphi_\ell}(Y_{(i)})(1-\delta_{(i)}) -\widehat{\gamma}_2^{\varphi_\ell}(Y_{(i)}),$$
$\widehat{\psi}=(\widehat{\psi}_{i1},\dots,\widehat{\psi}_{iL})^\top$ and $\widehat{\Psi}=(\widehat{\psi}_{1},\dots,\widehat{\psi}_{L})^\top$.
 Finally, the estimator of ${\Sigma}_\Psi$ is
 $\widehat{\Sigma}_{\Psi} =n^{-1}\sum_{i=1}^n \widehat{\Psi}\widehat\Psi^\top$ and that of $\Sigma$ is $\widehat{\Sigma}=\widehat W\widehat \Sigma_{\Psi}\widehat W^\top$.
 \section{Simulations} Let us evaluate the small sample performance of our estimation procedure in a simulation exercise. We consider the following data generating process (henceforth, DGP):
 \begin{align*}
 T&=0.5+X_2+X_3+U;\\
 X_2&=Z_2+V;\\
 U&=V+E,
 \end{align*}
 where $X_2, X_3, V$ and $E$ are mutually independent and distributed uniformly on the interval $[-1,1]$. In this DGP,
 $X=(1,X_2, X_3)^\top$, $Z=(1,Z_2,X_3)^\top$ and $\beta_0=(0.5,1,1)^\top$. The censoring variable $C$ is independent of the rest of the data and follows a unit exponential distribution. This results in a probability of censoring equal to $0.40$ (average over 1,000,000 replications). We want to evaluate the effect of the sample size $n$ on the accuracy of $\widehat{\beta}_2$. To do so, for $n=100, 1000,5000$, we compute the value of $\widehat{\beta}_2$ in $1,000$ Monte Carlo replications. In Table~\ref{samp}, we report the bias, the variance, the root mean-squared error (RMSE), the coverage of 95\% confidence intervals of  $\widehat{\beta}_2$, the width of these 95\% confidence intervals and the proportion of replications for which $\widehat{\beta}_2$ is significant. The estimator exhibits good performance even for low sample sizes.
\begin{table}[!h]
 
   \center {\small \begin{tabular}{|c|c|c|c|}
  \hline
 $n$ & $100$ &  $1,000$ &$5,000$ \\
  \hline
 Bias &-0.170& 0.035& 0.011 \\
Variance &0.134 & 0.014 &0.003\\
 RMSE &0.163 & 0.015 &0.003\\
  Coverage &0.88 & 0.89 &0.93\\
  Width &1.010 & 0.384&0.189\\
 $ \%$ Significant & 0.78 &1 &1\\
   \hline
      \end{tabular}    }
       \caption{ Bias, variance, RMSE, coverage of 95\% confidence intervals of  $\widehat{\beta}_2$, width of these 95\% confidence intervals and proportion of replications for which $\widehat{\beta}_2$ is significant; when $n=100, 1000, 5000$.}
             \label{samp}
\end{table}
In order to assess the effect of censoring on the estimator, we run additional simulations. The data generating process is the same with $n=1,000$ except that now $C=\rho+C_0$, where $C_0$ follows a $\text{Exp}(1)$. We assess the performance of the estimator when $n=1,000$ and $\rho =-1, -2, -3$ in $1,000$ Monte Carlo replications . These values of $\rho$ result in a probability of censoring of $0.61$ ($\rho =-1$), $0.80$ ($\rho =-2$), $0.91$ ($\rho=-3$) as computed through 1,000,000 simulations. In Table~\ref{cens}, we report the bias, the variance, the RMSE, the coverage of 95\% confidence intervals of  $\widehat{\beta}_2$, the width of these 95\% confidence intervals and the proportion of replications for which $\widehat{\beta}_2$ is significant. The estimator remains effective even when the proportion of censored observations is closed to $1$.

\begin{table}[!h]
 
   \center {\small \begin{tabular}{|c|c|c|c|}
  \hline
 $\rho$ & $-1$ &  $-2$ &$-3$ \\
  \hline
 Bias &-0.085& 0.127& 0.245\\
Variance &0.034 & 0.081 &0.290\\
 RMSE &0.041 & 0.097 &0.350\\
  Coverage &0.86 & 0.84 &0.83\\
  Width &0.56 & 0.784&1.20\\
 $ \%$ Significant & 0.98 &0.88 &0.71\\
   \hline
      \end{tabular}    }
       \caption{ Bias, variance, RMSE, coverage of 95\% confidence intervals of  $\widehat{\beta}_2$, width of these 95\% confidence intervals and proportion of replications for which $\widehat{\beta}_2$ is significant; when $\rho=-1,-2,-3$.}
             \label{cens}
\end{table}
 \section{Conclusion}
 This note shows how to adapt the usual 2SLS procedure to right censored outcomes. Our approach uses inverse probability of censoring weighting and retains the simplicity of the standard 2SLS estimator. Since, inverse probability of censoring weighting allows to estimate any moment, the same method could be applied to generalized method of moments estimators. We intend to develop this extension in future work.
 
\bibliographystyle{dcu}
\bibliography{ref}

\appendix
\section*{Appendix}
\subsection*{Proof of Theorem \ref{const}}
Notice that by assumption, $\E[ZZ^\top]$, $\E[ZX^\top]$ and $\E[ZY]=\E[ZX^\top]\beta_0$ are finite. Therefore, we can apply the the law of large numbers with censoring in \citet{Stute1993consistent}, which yields
 \begin{align*}
\sum_{i=1}^n  w_{(i)}Z_{(i)}Z_{(i)}^\top&\xrightarrow{\P}\E[ZZ^\top]\\
\sum_{i=1}^n  w_{(i)}Z_{(i)}X_{(i)}^\top&\xrightarrow{\P}\E[ZX^\top]\\
\sum_{i=1}^n  w_{(i)}Z_{(i)}Y_{(i)}&\xrightarrow{\P}\E[ZX^\top]\beta_0
 \end{align*}
 By the continuous mapping theorem and \eqref{closed-form}, we obtain 
 \begin{align*} \widehat{\beta}&\xrightarrow{\P}\left[\E[ZX^\top] \E[ZZ^\top] ^{-1}  \E[ZX^\top]\right]^{-1}\E[ZX^\top]\E[ZZ^\top]^{-1}  \E[ZX^\top]\beta_0=\beta_0.
 \end{align*}

\subsection*{Proof of Theorem \ref{AN}}

Given the proof of Theorem \ref{const}, it suffices to note that 
$$\left(\sum_{i=1}^n w_{(i)}Z_{(i)}Y_{(i)}\right) \xrightarrow{d} \mathcal{N}(0,\Sigma_{\Psi})$$
by Theorem 1.2 in \cite{stute1996distributional} and use Slutsky's theorem. Note that Theorem 1.2 in \cite{stute1996distributional} can be applied because Assumption \ref{censoring} (ii) implies that condition (1.3) in \cite{stute1996distributional} holds. 
\end{document}